# Comparisons of five indices for estimating local terrain surface roughness using LiDAR point clouds

Lei Fan
Department of Civil Engineering, Design School
Xi'an Jiaotong-Liverpool University
Suzhou, China
e-mail: Lei.Fan@xjtlu.edu.cn

*Abstract*—Terrain surface roughness is an abstract concept, and its quantitative description is often vague. As such, there are various roughness indices used in the literature, the selection of which is often challenging in applications. This study compared the terrain surface roughness maps quantified by five commonly used roughness indices, and explored their correlations for four terrain surfaces of distinct surface complexities. These surfaces were represented by digital elevation models (DEMs) constructed using airborne LiDAR (Light Detection and Ranging) data. The results of this study reveal the similarity in the global patterns of the local surface roughness maps derived, and the distinctions in their local patterns. The latter suggests the importance of considering multiple indices in the studies where local roughness values are the critical inputs to subsequent analyses.

## I. Introduction

Terrain surface roughness is an important index in geoscience to describe the complexity or variability of a terrain surface at a particular scale. It is used widely for many applications such as digital terrain analysis, Earth surface process simulations and terrain classifications [1]-[6]. Depending on the application requirements, global or local terrain roughness may be calculated. In the case where point cloud data are the source data, local terrain roughness is often calculated [7]-[11]. This is because such data have fine spatial resolution and record local terrain surface characteristics and the spatial changes in detail. Such data can readily be obtained using the measurement techniques such as Light Detection and Ranging (LiDAR) and Structure from Motion and Multi-View Stereo (SfM-MVS).

The definition of terrain surface roughness is often ambiguous. Roughness indices usually rely on a quantitative description of specific terrain characteristics changes, such as the degree of local undulation, the degree of local folds, or the degree of local abrupt changes. The indices commonly used include but are not limited to root mean square height (RMSH) [11], the standard deviation of residual elevation [9], the standard deviation of residual topography [9], the standard deviation of curvature [2], the standard deviation of slope [12]-[13], geostatistical analysis [14] and fractal dimension analysis [15].

Due to the wide range of applications of terrain surface roughness and the diversity of user needs, there is not yet a commonly accepted or preferred method for estimating it. In practice and research, people usually consider a commonly used terrain surface index at their preference or consider several indices to compare the results of interest for determining a more suitable index [7]-[11]. In the latter, the comparisons were generally based on a visual inspection [2], [9]. By now, there are few dedicated studies to investigate the quantitative correlations between terrain surface roughness indices, which motivated this study. In the study, five commonly used roughness indices were compared using four LiDAR point cloud datasets representing distinct terrain surfaces of varying surface complexities.

## II. Method

### A. Study data

Fig. 1 shows the study data, which consist of four airborne LIDAR point clouds representing the bare ground surfaces of a similar size of approximately 350 m by 350 m. These four terrain surfaces have different spatial variation characteristics: a hilly rough surface (an average data spacing of 0.67m), a hilly smooth surface (an average data spacing of 0.71 m), a flat rough surface (an average data spacing of 0.71 m), and a flat smooth surface (an average data spacing of 0.63m). These data were used in a previous study [16] and came from a large LiDAR dataset acquired by the National Airborne Laser Mapping Center in the USA at a volcanic field in central Nevada.

Each point cloud data has a global trend in elevation. To minimize its effect on the calculation of local surface roughness and to better visualize the spatial variations of the terrain surfaces, the global trend in each point cloud was removed by subtracting a best-fitting plane. This led to a set of residual elevations with a global mean elevation value of zero.

### B. DEM maps

Some surface roughness indices (e.g. RMSH) can be applied to either unstructured scattering point cloud data or structured/gridded data in the form of a digital elevation model (DEM). Meanwhile, there are also indices (e.g. standard deviation of slope) that are typically applicable to only gridded DEM maps. To keep the consistency of the comparisons

Accepted version of 2022 29th International Conference on Geoinformatics doi:10.1109/Geoinformatics57846.2022.9963877

This research was funded by XJTLU Key Program Special Fund (grant number KSF-E-40) and XJTLU Research Enhancement Funding (grant no. REF-21-01-003).

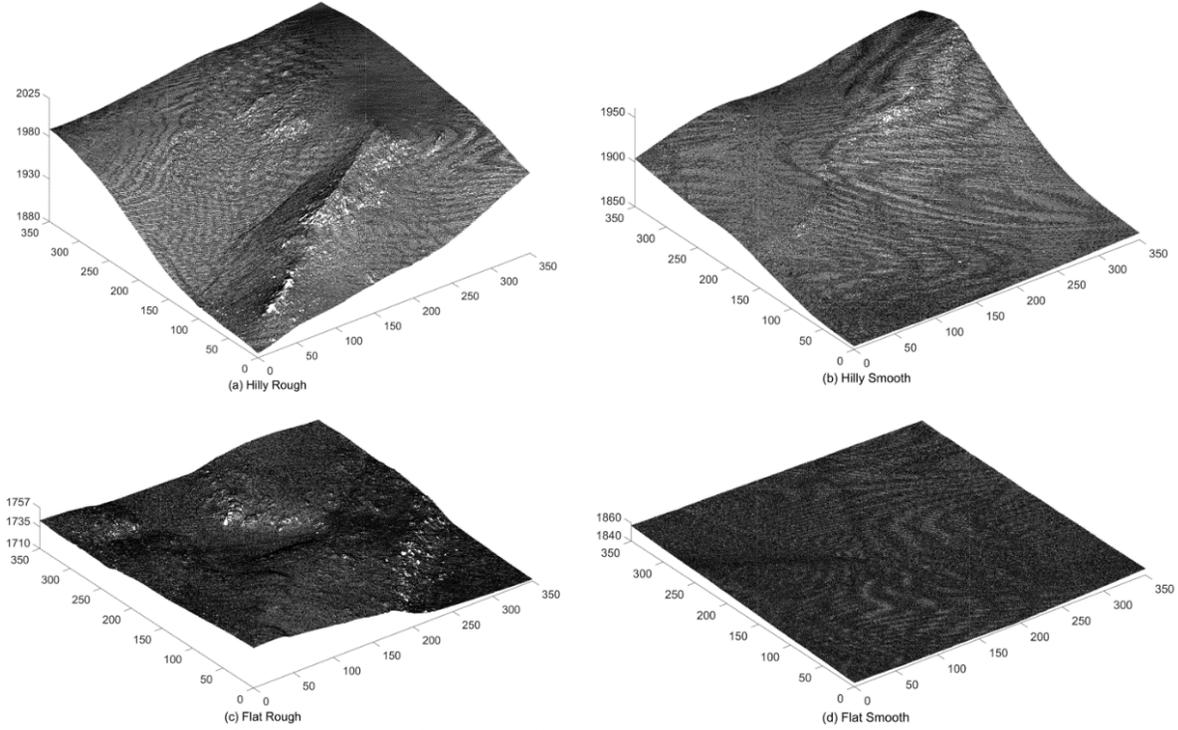

Fig. 1. The four point clouds considered: (a) Hilly Rough, (b) Hilly Smooth, (c) Flat Rough, (d) Flat Smooth.

between the roughness indices, gridded DEM maps were used to calculate local surface roughness in this study. The triangulation with a liner interpolation method [17] was used to produce those DEM maps of a spatial grid resolution of 1 m. The DEM maps constructed using the detrended point cloud data are shown in Fig. 2.

### C. Terrain surface roughness indices

The surface roughness indices considered in this study are elaborated in the following.

*1) RMSH*: It is one of the most commonly used indices for quantifying local surface roughness for scattered elevation data. The definition of RMSH is given in (1) [2], [11], [18].

$$\text{RMSE} = \sqrt{\frac{\sum_{i=1}^{n}(Z_i - \overline{Z})^2}{n-1}} \quad (1)$$

where $n$ represents the number of the data points selected; $Z_i$ is the elevation value of the $i^{\text{th}}$ data point; $\overline{Z}$ is the mean elevation value of all ($n$) data points selected.

*2) Standard deviation of locally detrended residual elevations ($\sigma_{\text{LDRE}}$)*: In this method, the local elevation data selected in a moving window are detrended linearly by fitting a best-fitting plane to obtain the residual elevation values of the data selected. The standard deviation of these residual elevations is calculated to represent the surface roughness.

*3) Standard deviation of residual topography ($\sigma_{\text{RT}}$):* The residual topography is the difference between the original topography (the original DEM) and the smoothed topography (the smoothed DEM) [9]. In this study, the elevation value at a grid location of the smoothed DEM was obtained by averaging the elevation values of its neighboring cells using a 5×5 moving window. As the original and the smoothed DEMs have the same spatial resolution, the residual topography is obtained by the arithmetic subtraction of the elevation values of the corresponding cells of the two DEMs. The standard deviation of the residual topography is calculated to represent the surface roughness.

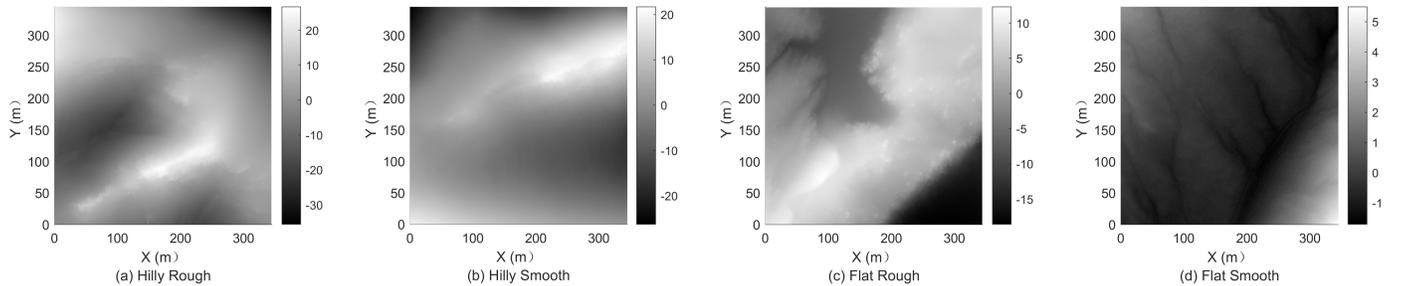

Fig. 2. The DEM maps of the detrended point cloud data: (a) Hilly Rough, (b) Hilly Smooth, (c) Flat Rough, (d) Flat Smooth.

*4) Standard deviation of slope ($\sigma_{slope}$):* To calculate the standard deviation of slope, it is necessary to first calculate the slope values. Slope is defined as the rate of change of terrain elevations, which is expressed in (2).

$$\text{slope} = \tan^{-1}\left(\sqrt{\left(\frac{dz}{dx}\right)^2 + \left(\frac{dz}{dy}\right)^2}\right) \quad (2)$$

where $dz/dx$ and $dz/dy$ represent the rate of change in the $x$ and the $y$ directions, respectively, for the cell of interest.

Based on a DEM, slope is often calculated using the 3×3 moving window shown in Fig. 3. $dz/dx$ and $dz/dy$ for the central cell ($Z_5$) are calculated using (3) and (4), respectively, in which $L$ represents the cell size. When there are cells in the neighborhood (e.g. at the edge of a DEM) that do not contain elevation data, those cells are assumed to have the same elevation value as the central cell. This is useful for the cells at the edge of a DEM raster to ensure that the slope map has the same spatial extent as the DEM map.

$$dz/dx = [(Z_3+2Z_6+Z_9) - (Z_1+2Z_4+Z_7)]/(8L) \quad (3)$$

$$dz/dy = [(Z_7+2Z_8+Z_9) - (Z_1+2Z_2+Z_3)]/(8L) \quad (4)$$

*5) Standard deviation of curvature ($\sigma_{curvature}$):* Curvature is calculated by computing the second derivative of a DEM raster map. The moving window used is the same as that (i.e. Fig. 3) for calculating slope. There are various methods for calculating curvature. In this study, the one proposed by Zevenbergen and Thorne [19]-[20] was adopted, which is described in (5) - (7). Similar to the calculation of slope, the non-value cells in the neighborhood are assumed to have the same elevation value as the central cell. Once the curvature raster map is obtained, the standard deviation of curvature is calculated to represent the surface roughness.

$$\text{curvature} = 2E + 2D \quad (5)$$

where $E$ and $D$ are given in (6) and (7), respectively, using the elevation values within the moving window shown in Fig. 3.

$$D = [(Z_4 + Z_6)/2 - Z_5]/L^2 \quad (6)$$

$$E = [(Z_2 + Z_8)/2 - Z_5]/L^2 \quad (7)$$

### D. The spatial scale for local terrain surface roughness

The local terrain surface roughness is often calculated at a particular scale defined by users, which is implemented using non-overlapping moving windows. Following a survey of the typical scales used in the literature for computing local surface roughness using DEM maps derived from point cloud data, the following set of non-overlapping moving window sizes (3×3, 5×5, 7×7, 9×9 and 11×11) were considered in this study to look into effects of the spatial scale on the local surface roughness maps produced by the different roughness indices.

| $Z_1$ | $Z_2$ | $Z_3$ |
|---|---|---|
| $Z_4$ | $Z_5$ | $Z_6$ |
| $Z_7$ | $Z_8$ | $Z_9$ |

Fig. 3. The 3×3 moving window for calculating slopes using a DEM

### III. RESULTS

Fig. 4 shows the local terrain surface roughness maps of the five surface roughness indices considered, based on a spatial scale of 5×5 cells as an example. This is also a typical spatial scale that was commonly considered in the literature. Those derived using the other spatial scales are not shown here due to the page limit. In Fig. 4, each individual column of the plots represents a different roughness index while each row represents one of the four terrain surfaces considered. It should be noted that the roughness values shown in Fig. 4 were normalized to range from 0 to 1 for a better visual comparison.

It was observed that most of the roughness indices produced local terrain surface roughness maps of similar global patterns, justifying the use of these indices in the literature. To confirm this, the correlation coefficient $r$ for each pair of roughness maps/images compared was calculated using (8) for the Hilly Rough data as an example. The values of correlation coefficient are shown in TABLE I where only those in the upper triangular are shown due to due to symmetry. Apart from RMSH, large correlations were found between two roughness maps derived using the other indices. The largest correlation value (i.e. 0.95) was observed for the pair $\sigma_{RT}$ and $\sigma_{curvature}$, which is consistent with the global patterns of the roughness maps for these two indices as shown in Fig. 4.

$$r = \frac{\sum_m \sum_n (A_{mn} - \overline{A})(B_{mn} - \overline{B})}{\sqrt{\left(\sum_m \sum_n (A_{mn} - \overline{A})^2\right)\left(\sum_m \sum_n (B_{mn} - \overline{B})^2\right)}} \quad (8)$$

where $A$ and $B$ represent the values of the roughness images compared; the subscript $m$ and $n$ refer to the pixel location in the images; $\overline{A}$ and $\overline{B}$ represent the mean value.

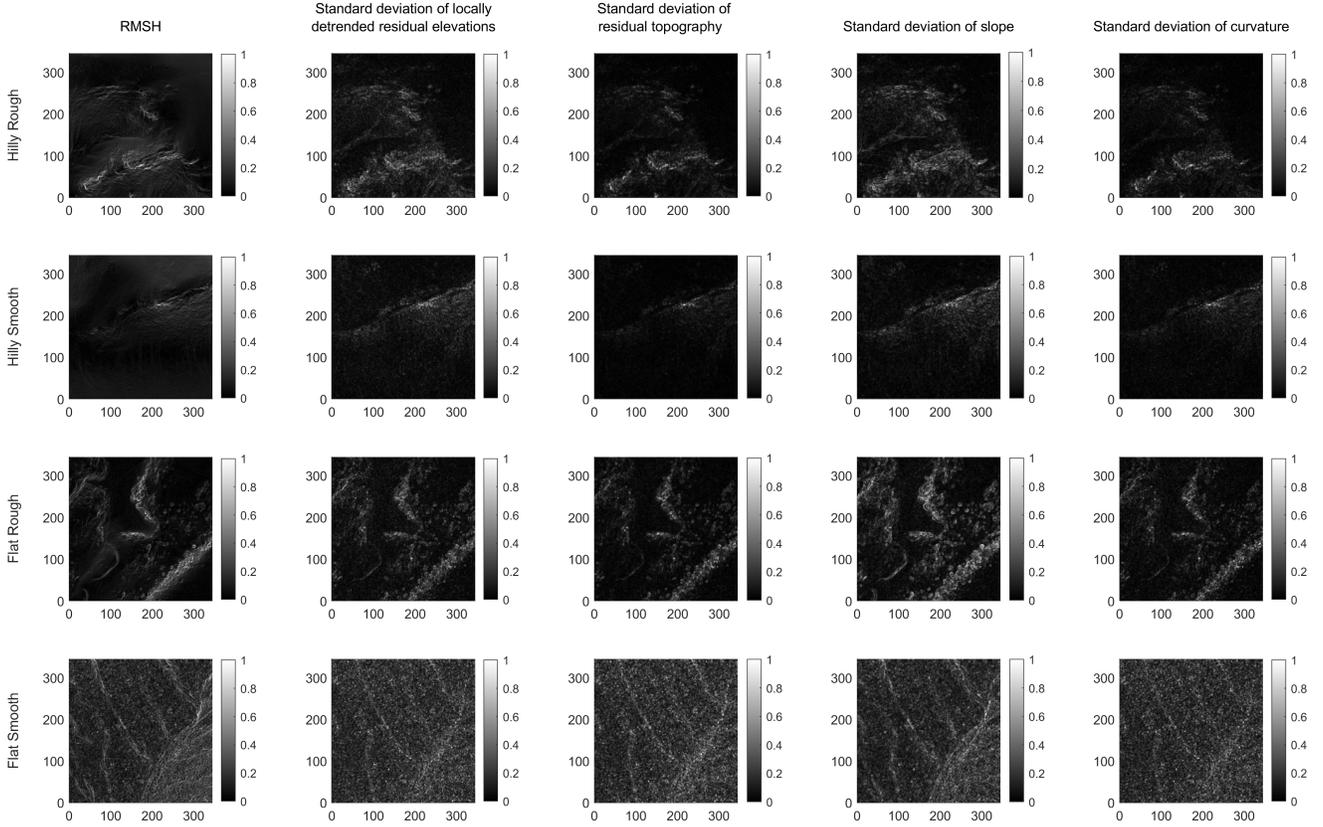

Fig. 4. The local roughness maps for the four terrain surfaces considered.

TABLE I. THE CORRELATION COEFFICIENT FOR EACH PAIR OF ROUGHNESS MAPS COMPARED, BASED ON THE HILLY ROUGH DATA

| Roughness Index | Roughness Index | | | | |
|---|---|---|---|---|---|
| | RMSH | $\sigma_{LDRE}$ | $\sigma_{RT}$ | $\sigma_{slope}$ | $\sigma_{curvature}$ |
| RMSH | 1 | 0.59 | 0.72 | 0.52 | 0.71 |
| $\sigma_{LDRE}$ | | 1 | 0.83 | 0.81 | 0.88 |
| $\sigma_{RT}$ | | | 1 | 0.76 | 0.95 |
| $\sigma_{slope}$ | | | | 1 | 0.74 |
| $\sigma_{curvature}$ | | | | | 1 |

However, the distributions of the roughness values locally vary between those indices. Some (e.g. $\sigma_{RT}$ and $\sigma_{curvature}$) exhibit similar local distributions while some others (e.g. RMSH and $\sigma_{slope}$) are notably different. These variations between the indices suggest that the choice of roughness indices can affect the results of a subsequent analysis because many quantitative studies involving local surface roughness are often based on local roughness values. As such, it would be useful to consider a sensitivity analysis using multiple roughness indices in the quantitative studies where the local roughness values are the critical inputs.

Another means of quantitative estimation of the correlation between roughness maps is the coefficient of determination $R^2$, which was calculated using all pixel values of one roughness map against the corresponding pixel values of another roughness map being compared. Two representative datasets (i.e. the most complex surface Hilly Rough and the least complex surface Flat Smooth) were used for the analysis, in which five local scales (i.e. moving windows of 3×3, 5×5, 7×7, 9×9 and 11×11 cells) were used. The results are summarized graphically in Fig. 5. The overall correlations between two indices compared appeared to be slightly larger for the rougher surface (i.e. Hilly Rough).

The correlation between RMSH and any one of the other four roughness indices was found to be relatively small, as illustrated in Fig. 5. This can also be inferred from the roughness maps in Fig. 4. It is interesting to observe in Fig. 4 that the roughness values of RMSH are more spatially coherent than those of the other indices. The likely reason for these is the presence of local elevation trends, which can lead to stronger local spatial autocorrelation in the maps of RMSH. This explanation can partly be justified by the map of $\sigma_{LDRE}$, which share the same algorithm as RMSH but use locally de-trended elevation values.

The effects of the spatial scales on the correlations between the roughness indices are shown in Fig. 5. A mixed behavior was observed. For the pairs of roughness indices compared, a larger window size may increase or reduce their correlations. However, based on the results shown in Fig. 5, it appears that the effect of the window size is smaller for the rougher terrain surface (i.e. Hilly Rough).

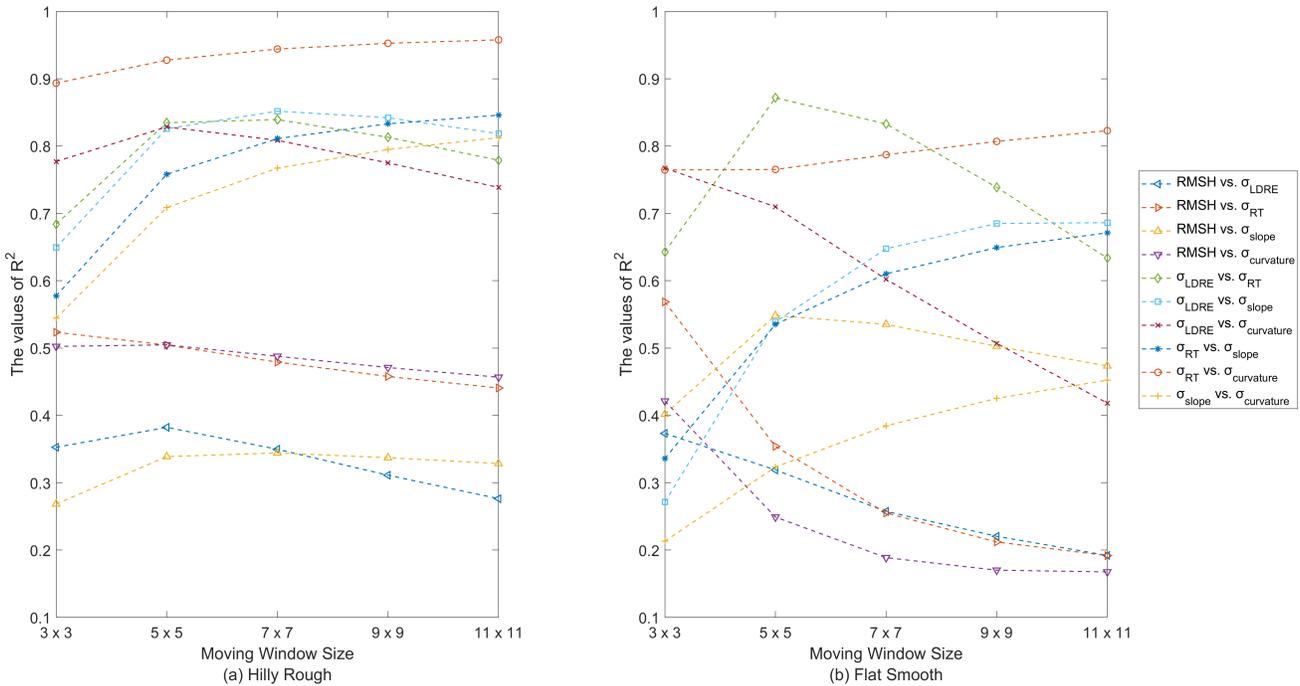

Fig. 5. The coefficient of determination for paired roughness indices under different window sizes: (a) Hilly Rough and (b) Flat Smooth.

## IV. Conclusions

For the four terrain surfaces of varying complexities, the local roughness maps produced by the roughness indices considered showed similar global patterns, demonstrating their capability of characterizing local terrain surface roughness and their suitability for applications. However, it was also found that the local distributions of the roughness values varied notably between some roughness indices, suggesting the importance of considering multiple roughness indices in the studies where local roughness values are the input for a subsequent analysis. This is particularly the case for RMSH as its correlations with the other four indices considered were found to be small. The effects of local scales (i.e. moving window sizes) on the correlation between two roughness indices compared were briefly investigated, which showed a mixed behavior (i.e. either decreased or increased the correlation).

## Acknowledgment

LiDAR data access is based on [LiDAR, ground] services provided by the OpenTopography Facility with support from the National Science Foundation under NSF Award Numbers 1226353 & 1225810. Lidar data acquisition completed by the National Center for Airborne Laser Mapping (NCALM - http://www.ncalm. org). NCALM funding provided by NSF's Division of Earth Sciences, Instrumentation and Facilities Program. EAR-1043051. https://doi.org/10.5069/G9PR7SX0